\def\timestamp{%
Time-stamp: <universal-homeomorphism.tex: Sunday 05-09-2021 at 18:49:50 (cest)>}
\def\stripname Time-stamp: <#1 #2>{#2}
\edef\filedate{\expandafter\stripname\timestamp}

\documentclass{amsart}
\usepackage{amsrefs}

\newcommand\Omegaseq[2][1]{\langle {#2}_\alpha:\alpha\in\omega_{#1}\rangle}

\newcommand\Aut{\mathsf{Aut}}

\newcommand\cl{\operatorname{cl}}
\newcommand\St{\operatorname{St}}
\newcommand\fin{\mathrm{fin}}
\newcommand\calC{\mathcal{C}}
\newcommand\calU{\mathcal{U}}

\newcommand\cee{\mathfrak{c}}

\newcommand\axiom{\mathsf}
\newcommand\CH{\axiom{CH}}

\DeclareMathSymbol\B0{AMSb}{`B}
\DeclareMathSymbol\N0{AMSb}{`N}
\DeclareMathSymbol\Z0{AMSb}{`Z}
\newcommand\Nstar{\N^*}

\DeclareMathSymbol\restr\mathbin{AMSa}{"16}
\DeclareMathSymbol\le    \mathrel{AMSa}{"36}
\DeclareMathSymbol\ge    \mathrel{AMSa}{"3E}

\begin{document}

\title{Universal autohomeomorphisms of $\Nstar$}

\author{Klaas Pieter Hart}

\address{Faculty EEMCS\\TU Delft\\
         Postbus 5031\\2600~GA {} Delft\\the Netherlands}
\email{k.p.hart@tudelft.nl}
\urladdr{http://fa.ewi.tudelft.nl/\~{}hart}

\author{Jan van Mill}
\address{KdV Institute for Mathematics\\
         University of Amsterdam\\
         P.O. Box 94248\\
         1090~GE {} Amsterdam\\
         The Netherlands}
\email{j.vanmill@uva.nl}
\date{\filedate}


\begin{abstract}
We study the existence of universal autohomeomorphisms of $\Nstar$.
We prove that $\CH$ implies there is such an autohomeomorphism and show
that there are none in any model where all autohomeomorphisms
of~$\Nstar$ are trivial.  
\end{abstract}

\dedicatory{To the memory of Cor Baayen, 
            who taught us many things}

\maketitle

\section*{Introduction}

This paper is concerned with universal autohomeomorphisms on~$\Nstar$, 
the \v{C}ech-Stone remainder of~$\N$.

In very general terms we say that an autohomeomorphism~$h$ on a space~$X$
is universal for a class of pairs $(Y,g)$, where $Y$~is a space and $g$~is
an autohomeomorphism of~$Y$, if for every such pair there
is an embedding $e:Y\to X$ such that $f\circ e=e\circ g$, that is,
$h$~extends the copy of~$g$ on~$e[Y]$.
 
In \cite{MR0172826}*{Section~3.4} one finds a general way of finding
universal autohomeomorphisms.
If $X$~is homeomorphic $X^\omega$ then the shift 
mapping $\sigma:X^\Z\to X^\Z$ defines a universal autohomeomorphism for
the class of all pairs $(Y,g)$, where $Y$~is a subspace of~$X$.
One embeds $Y$ into~$X^\Z$ by mapping each~$y\in Y$ to the 
sequence $\langle g^n(y):n\in\Z\rangle$. 

Thus, the Hilbert cube carries an autohomeomorphism
that is universal for all autohomeomorphisms of separable metrizable
spaces and the Cantor set carries one for all autohomeomorphisms of 
zero-dimensional separable metrizable spaces.
Likewise the Tychonoff cube~$[0,1]^\kappa$ carries an autohomeomorphism
that is universal for all autohomeomorphisms of completely regular spaces
of weight at most~$\kappa$, and the Cantor cube~$2^\kappa$ has a universal 
autohomeomorphism for all zero-dimensional such spaces.   

Our goal is to have an autohomeomorphism~$h$ on~$\Nstar$ that is universal
for all autohomeomorphisms of all \emph{closed} subspaces of~$\Nstar$.
The first result of this paper is that there is no trivial universal
autohomeomorphism of~$\Nstar$, and hence no universal autohomeomorphism
at all in any model where all autohomeomorphisms of~$\Nstar$ are trivial.
On the other hand, the Continuum Hypothesis implies that there is
a universal autohomeomorphism of~$\Nstar$.
The proof of this will have to be different from the results mentioned 
above because $\Nstar$ is definitely not homeomorphic to its 
power~$(\Nstar)^\omega$; it will use group actions and a homeomorphism 
extension theorem.

We should mention the dual notion of universality where one requires
the existence of a surjection $s:X\to Y$ such that $g\circ s=s\circ h$.
For the space~$\Nstar$ this was investigated thoroughly in~\cite{MR4013982}
for general group actions.

\section{Some preliminaries}

Our notation is standard.
For background information on~$\Nstar$ we refer to~\cite{MR776630}.

We denote by~$\Aut$ the autohomeomorphism group of~$\Nstar$.
We call a member~$h$ of~$\Aut$ \emph{trivial} if there are cofinite
subsets~$A$ and~$B$ of~$\N$ and a bijection~$b:A\to B$ such that
$h$~is the restriction of~$\beta b$ to~$\Nstar$.

In both sections we shall use the $G_\delta$-topology on a given 
space~$(X,\tau)$; this is the topology~$\tau_\delta$ on~$X$ generated by
the family of all $G_\delta$-subsets in the given space.
It is well-known that $w(X,\tau_\delta)\le w(X,\tau)^{\aleph_0}$; we shall
need this estimate in Section~\ref{sec.CH}.

\section{What if all autohomeomorphisms are trivial?}

To begin we observe that fixed-point sets of trivial autohomeomorphism 
of~$\Nstar$ are clopen.
Therefore, to show that no trivial autohomeomorphism is universal
it would suffice to construct a compact space that can be embedded
into~$\Nstar$ and that has an autohomeomorphism whose fixed-point set
is not clopen.

\subsection*{The example}
We let $L$ be the ordinal~$\omega_1+1$ endowed with its $G_\delta$-topology.
Thus all points other than~$\omega_1$ are isolated and the neighbourhoods
of~$\omega_1$ are exactly the co-countable sets that contain it.
Then $L$~is a $P$-space of weight~$\aleph_1$ and hence, 
by the methods in~\cite{MR674103}*{Section~2},
its \v{C}ech-Stone compactification~$\beta L$ can be embedded into~$\Nstar$.

We define $f:L\to L$ such that $\omega_1$ is the only fixed point 
of~$\beta f$.
We put
\begin{align*}
  f(\omega_1)&=\omega_1,\\
  f(2\cdot\alpha)&=2\cdot\alpha+1, \text{ and }\\
  f(2\cdot\alpha+1)&=2\cdot\alpha.
\end{align*}
This defines a continuous involution on~$L$.

If $p\in\beta L\setminus L$ then $p\in\cl\alpha$ for some~$\alpha<\omega_1$
and then either $E=\{2\cdot\beta:\beta<\alpha\}$ or
$O=\{2\cdot\beta+1:\beta<\alpha\}$ belongs to the ultrafilter~$p$.
But $f[E]\cap E=\emptyset=f[O]\cap O$, hence $\beta f(p)\neq p$. 

Since $\omega_1$ is not an isolated point of~$\beta L$, no matter how
this space is embedded into~$\Nstar$ there is no trivial autohomeomorphism
of~$\Nstar$ that would extend~$\beta f$.

\section{The Continuum Hypothesis}\label{sec.CH}

Under the Continuum Hypothesis the space $\Nstar$ is generally very 
well-behaved and one would expect it to have a universal autohomeomorphism 
as well.
We shall prove that this is indeed the case.
We need some well-known facts about closed subspaces of~$\Nstar$.

First we have Theorem~1.4.4 from~\cite{MR776630} which characterizes 
the closed subspaces of~$\Nstar$ under~$\CH$:
they are the compact zero-dimensional $F$-spaces of weight~$\cee$, and,
in addition: every closed subset of~$\Nstar$ can be re-embedded as a 
nowhere dense closed $P$-set.

Second we have the homeomorphism extension theorem from~\cite{MR1277871}:
$\CH$~implies that every homeomorphism between nowhere dense closed $P$-sets
of~$\Nstar$ can be extended to an autohomeomorphism of~$\Nstar$.

\subsection*{Step 1}

We consider the natural action of $\Aut$ on~$\Nstar$, that is
the map $\sigma:\Aut\times\Nstar\to\Nstar$ given by $\sigma(f,p)=f(p)$.
This action is continuous when $\Aut$~carries the compact-open topology~$\tau$
and hence also when $\Aut$~carries the $G_\delta$-modification~$\tau_\delta$ 
of~$\tau$.
For the rest of the construction we consider the topology~$\tau_\delta$.

Using this action we define an autohomeomorphism 
$h:\Aut\times\Nstar\to\Aut\times\Nstar$
by $h(f,p)=(f,f(p))$.
The map~$h$ is continuous because its two coordinates are and it is 
a homeomorphism because its inverse $(f,p)\mapsto(f,f^{-1}(p))$ is continuous as
well. 

Now if $X$ is a closed subset of~$\Nstar$ and $g:X\to X$ is an 
autohomeomorphism then we can re-embed $X$ as a nowhere dense closed $P$-set
and we can then find an $f\in\Aut$ such that $f\restr X=g$. 
We transfer this embedded copy of~$X$ to~$\{f\}\times\Nstar$ 
in~$\Aut\times\Nstar$; for this copy of~$X$ we then have $h\restr X=g$.
It follows that $h$~satisfies the universality condition.

\subsection*{Step 2}
We embed $\Aut\times\Nstar$ into~$\Nstar$ in such a way that
there is an autohomeomorphism~$H$ of~$\Nstar$ such that
$H\restr(\Aut\times\Nstar)=h$.
Then $H$~is the desired universal autohomeomorphism of~$\Nstar$. 

\smallskip
To this end we list a few properties of this product.

\subsubsection*{Weight}
The weight of the product is equal to~$\cee$, as both factors have 
weight~$\cee$.
For $\Nstar$ this is clear and for $\Aut$ this follows because the
topology~$\tau$ has weight~$\cee$ and one obtains a base for~$\tau_\delta$
by taking the intersections of all countable subfamilies of a base for~$\tau$.

\subsubsection*{Zero-dimensional and $F$}
The product is a zero-dimensional $F$-space as the product of the 
$P$-space~$\Aut$ and the compact zero-dimensional $F$-space~$\Nstar$, 
see~\cite{MR234407}*{Theorem~6.1}.

\subsubsection*{Strongly zero-dimensional}

The product $\Aut\times\Nstar$ is not compact, but we shall construct a 
compactification of it that is also a zero-dimensional $F$-space of 
weight~$\cee$.

For this we need to prove that $\Aut\times\Nstar$ is actually strongly 
zero-dimensional.
We prove more: the product is ultraparacompact, meaning that every open cover 
has a \emph{pairwise disjoint} open refinement.

\smallskip
Let $\calU$ be an open cover of the product consisting of basic clopen 
rectangles.

For each~$f\in\Aut$ there is a finite subfamily $\calU_f$ of~$\calU$
that covers $\{f\}\times\Nstar$, say 
$\calU_f=\{C_i\times D_i:i<k_f\}$.
Let $C_f=\bigcap_{i<k}C_i$ and $D_{f,i}=D_i\setminus\bigcup_{j<i}D_j$ for $i<k_f$.
Then $\calC_f=\{C_f\times D_{f,i}:i<k_f\}$ is a disjoint family of clopen 
rectangles that covers~$\{f\}\times\Nstar$ and refines~$\calU$.

Because $\Aut$ has weight~$\cee$, and we assume~$\CH$, there is a sequence
$\Omegaseq f$ in~$\Aut$ such that $\{C_{f_\alpha}:\alpha\in\omega_1\}$ 
covers~$\Aut$. 
Next we let $V_\alpha=C_{f_\alpha}\setminus\bigcup_{\beta<\alpha}C_{f_\beta}$ for 
all~$\alpha$.
Because $\Aut$~is a $P$-space the family $\{V_\alpha:\alpha\in\omega_1\}$ is
a disjoint open cover of~$\Aut$.
 
The family $\{V_\alpha\times D_{f_\alpha,i}:i<k_{f_\alpha}, \alpha\in\omega_1\}$ then
is a disjoint open refinement of~$\calU$.

\subsubsection*{A compactification}
To complete Step~2 we construct a compactification of~$\Aut\times\Nstar$
that is a zero-dimensional $F$-space of weight~$\cee$ and that has an 
autohomeomorphism that extends~$h$.
The \v{C}ech-Stone compactification would be the obvious canditate, were
it not for the fact that its weight is equal to~$2^\cee$.
More precisely, using some continuous onto function from $(\Aut,\tau)$
onto~$[0,1]$ one obtains a clopen partition of $(\Aut,\tau_\delta)$ of
cardinality~$\cee$. 
This shows that $\beta(\Aut\times\Nstar)$ admits a continuous surjection
onto the space~$\beta\cee$ (where $\cee$~carries the discrete topology).

To create the desired compactification we build, either by transfinite
recursion or by an application of the L\"owenheim-Skolem theorem,
a subalgebra~$\B$ of the algebra of clopen subsets of~$\Aut\times\Nstar$
that is closed under~$h$ and~$h^{-1}$, of cardinality~$\cee$,
and that has the property that 
for every pair of countable subsets~$A$ and~$B$ of~$\B$ such 
that $a\cap b=\emptyset$ whenever $a\in A$ and $b\in B$ 
there is a $c\in\B$ such that $a\subseteq c$
and $c\cap b=\emptyset$ for all $a\in A$ and $b\in B$.
The latter condition can be fulfilled because $\Aut\times\Nstar$~is 
an $F$-space --- $\bigcup A$ and~$\bigcup B$ have disjoint closures ---
and strongly zero-dimensional --- the closures can be separated using
a clopen set. 

The Stone space $\St(\B)$ of~$\B$ is then a compactification 
of~$\Aut\times\Nstar$ that is a compact zero-dimensional $F$-space of
weight~$\cee$, with an autohomeomorphism~$\bar h$ that extends~$h$.
We embed $\St(\B)$ into~$\Nstar$ as a nowhere dense $P$-set and 
extend~$\bar h$ to an autohomeomorphism~$H$ of~$\Nstar$.

\end{document}